\documentclass[11pt]{article}
\usepackage{amsthm}
\usepackage{amsmath}
\usepackage{amssymb}
\usepackage{amsfonts,graphicx}
\usepackage[all]{xy}
\usepackage{bezier}
\numberwithin{equation}{section}

\newtheorem{defn}{Definition}[section]
\newtheorem{thm}{Theorem}[section]
\newtheorem{prop}{Proposition}[section]
\newtheorem{lem}{Lemma}[section]
\newtheorem{rmk}{Remark}[section]

\theoremstyle{plain}

\begin{document}
\title{The moduli space of  holomorphic chains of rank one over a compact Riemann surface}
\author{JinHyung To}
\maketitle{}
\begin{abstract}
A holomorphic chain on a compact Riemann surface is a tuple of vector bundles together with homomorphisms between them. We show that the moduli space of holomorphic chains of rank one is identified with a fiber product of projective space bundles. We compute the Euler characteristic of the moduli space. The stability of chains involves real vector parameters. We also show that the variation of parameters corresponds to the characters of $G_m$.
\end{abstract}

%\tableofcontents

% Section 1 %%%%%%%%%%%%%%%%%%%%%%%%%%%%%%%%%%%%%%%%%%%

\section{Introduction}

% category of chains

Let $X$ be a compact Riemann surface of genus $g\geq 0$.
\begin{defn} A holomorphic chain $C$ on $X$ is a tuple
$(E_0,...,E_n;\phi_1,...,\phi_n)$ of vector bundles together with the homomorphisms $\phi_i:E_i\rightarrow E_{i-1}$ for $i=1,...,n$. The
tuple $\boldsymbol{t}=(rk(E_0),...,rk(E_n);deg(E_0),...,deg(E_n))$
is called the type of the holomorphic chain $C$. A subchain of  $C$ is a holomorphic
chain $(F_0,...,F_n;\phi_1|F_1,...,\phi_n|F_n),$ where $F_i\subseteq
E_i$ is a subsheaf such that $\phi_i(F_i)\subseteq F_{i-1}$ for
$1\leq i\leq n$.
\end{defn}

\begin{rmk} Since $E_i$ is a torsion free sheaf over a compact Riemann surface, the subsheaf
$F_i$ is torsion free and it is a vector bundle.
\end{rmk}

We call a holomorphic chain just a chain to be brief. There is a concept of
stability for a chain involving a real vector parameter. Given a chain $C$ of type
$(r_0,...,r_n;d_0,...,d_n)$ and a vector $\boldsymbol{\alpha} =
(\alpha_0,\alpha_1,...,\alpha_n) \in \mathbb{R}^{n+1}$, the
$\boldsymbol{\alpha}$-slope of $C$ is the fraction
\[
\mu_{\boldsymbol{\alpha}}(C)=\frac{\sum_{i=0}^n(d_i + \alpha_i
r_i)}{\sum_{i=0}^n r_i}.
\]
The chain $C$ is
$\boldsymbol{\alpha}$-stable ($\boldsymbol{\alpha}$-semistable) if
\[
\mu_{\boldsymbol{\alpha}}(C') < \mu_{\boldsymbol{\alpha}}(C) \, (\mu_{\boldsymbol{\alpha}}(C') \le \mu_{\boldsymbol{\alpha}}(C))
\]
for all proper subchains $C'$, i.e., $C' \ne (0,\cdots,0;0,\cdots,0)$ and $C' \ne C$.

\begin{rmk}
i) A subsheaf $F$ of a vector bundle is contained in a subbundle $\bar{F}$ with $deg(F) \leq
deg(\bar{F})$. Thus semistability can be checked against all subchains
composed of subbundles.

\noindent ii) Let $\boldsymbol{\alpha}=(\alpha_0,...,\alpha_n)$. Set
$\boldsymbol{\beta} = (\alpha_0 + \alpha,...,\alpha_n + \alpha)$. Then
\[
\mu_{\boldsymbol{\beta}}(C) = \frac{\sum_{i = 0}^n(d_i + (\alpha_i + \alpha)r_i)}{\sum_{i = 0}^n r_i}
=\mu_{\boldsymbol{\alpha}}(C) + \boldsymbol{\alpha};
\]
so translating by a vector $(\alpha,...,\alpha)$ does not change $\boldsymbol{\alpha}$-stability, i.e., the chain $C$ is $\boldsymbol{\alpha}$-(semi)stable if and only if it is $\boldsymbol{\beta}$-(semi)stable. So we may assume that $\alpha_0 =
0$ by setting $\alpha = -\alpha_0$. Thus the $\boldsymbol\alpha$-slope is reduced to
\[
\mu_{\boldsymbol{\alpha}}(C)=\frac{\sum_{i=0}^n d_i + \sum_{i=1}^n
\alpha_i r_i}{\sum_{i=0}^nr_i}
\]
for $\boldsymbol{\alpha} \in \mathbb{R}^n$. From now on we use this reduced $\boldsymbol\alpha$-slope.
\end{rmk}

Holomorphic chains form the category of representations of quivers $\bullet
\rightarrow \bullet \rightarrow ... \rightarrow \bullet$ in the 
category of vector bundles on $X$. A quiver is an oriented graph
consisting of the set of vertices and the set of arrows. A
representation in the category of vector bundles means the
assignments of the vector bundles to each vertex and of the maps to each
arrow. These are called \emph{quiver bundles}.

A map from $C = (E_0,...,E_n;\phi_1,...,\phi_n)$ to $C ' =
(F_0,...,F_n;\psi_1,...,\psi_n)$ is a tuple of vertical maps $f_i :
E_i \rightarrow F_i$ for $i = 0,..., n$ such that $\psi_i f_i = f_{i-1} \phi_i$ for $i
= 1,...,n$. It is an isomorphism if all $f_i$'s are isomorphisms.

A holomorphic chain $Q = (Q_0,...,Q_n;\tilde{\phi_1},...,\tilde{\phi_n})$ 
is called a \emph{quotient chain} of \\ $C = (E_0,...,E_n;\phi_1,...,\phi_n)$ if
there exists a surjective homomorphism $\pi_i:E_i\rightarrow Q_i$
for $i=0,...,n$ such that $\tilde{\phi_i}\pi_i=\pi_{i-1}\phi_i$ for $i=1,...,n$. If $C$ is $\boldsymbol\alpha $-stable, then by the additivity of degree of vector bundles, $\mu_{\boldsymbol\alpha}(C) < \mu_{\boldsymbol\alpha}(Q)$.
The category of holomorphic chains does not form an abelian category 
as the category of vector bundles does not form an abelian category.
Like vector bundles, however, 
the subcategory of $\boldsymbol{\alpha}$-semistable chains with the same slope
$\mu_{\boldsymbol{\alpha}}$ form an abelian category. Moreover, since it is
artinian and noetherian, the Jordan-H\"{o}lder Theorem holds in the category of chains. The simple objects are precisely the $\boldsymbol{\alpha}$-stable chains.

Let $C$ be an $\boldsymbol{\alpha}$-semistable chain of slope $\mu_{\boldsymbol{\alpha}}$. If $0 =C_0 \subsetneq C_1 \subsetneq ...\subsetneq C_m = C $ is a Jordan-H\"{o}lder filtration, i.e., $C_i/C_{i-1}$ is $\boldsymbol{\alpha}$-stable with slope $\mu_{\boldsymbol{\alpha}}(C_i/C_{i-1})=\mu_{\boldsymbol{\alpha}}$
for $i = 1,...,m$, then we define
\[
\text{gr}(C):=\bigoplus_{i=1}^m(C_i/C_{i-1})
\] called the
\emph{graduation} of $C$. Then the two $\boldsymbol{\alpha}$-semistable
chains $C$ and $C'$ are called \emph{S-equivalent} if $\text{gr}(C) \cong
\text{gr}(C')$. In particular, gr$(C)=C$ if $C$ is $\boldsymbol{\alpha}$-stable.

% moduli spaces of chains

Finally, the moduli spaces of holomorphic chains are constructed in \cite{Sch}. 
Indeed the moduli spaces of representations of more general quivers in the category of sheaves on a projective manifold are constructed. They are projective schemes and hence the moduli spaces of 
holomorphic chains are also projective schemes.

\noindent\textbf{Acknowledgements.} I would like to thank my advisor Steven Bradlow for proposing
this question and for valuable discussions.

% Section 2. %%%%%%%%%%%%%%%%%%%%%%%%%%%%%%%%%%%%%%%%%%%

\section{The standard chamber of the $\boldsymbol{\alpha}$ vector parameters}

A chain $C = (E_0, E_1, ..., E_n; \phi_1, ..., \phi_n)$ has a subchain
\[
C_i = (E_0, E_1, ..., E_i, 0, ..., 0;\phi_1, ..., \phi_i,0,...,0)
\]
called the i-th standard subchain for $i = 0, ..., n-1$. 
Given a vector parameter $\boldsymbol{\alpha} = (\alpha_1, ..., \alpha_n)$ and 
a type $\boldsymbol{t} = (r_0, r_1, ..., r_n; d_0, d_1, ..., d_n)$ of $C$, if $C$ is $\boldsymbol{\alpha}$-stable, then
$\mu_{\boldsymbol{\alpha}}(C_i) < \mu_{\boldsymbol{\alpha}}(C)$ for all $i = 0,...,n-1$. 
These imply the following $n$ inequalities:
\[
\sum_{k=i+1}^nr_k \sum_{k=1}^i \alpha_k r_k - \sum_{k=0}^i r_k  \sum_{k=i+1}^n \alpha_k r_k  < 
\sum_{k=0}^i r_k \sum_{k=i+1}^n d_k  - \sum_{k=i+1}^n r_k \sum_{k=0}^i d_k 
\]
for $ i = 0, ..., n-1$. Denote by $h_{i+1}$ the hyperplane in $\mathbb{R}^n$ determined by $\mu_{\boldsymbol{\alpha}}(C_i) = \mu_{\boldsymbol{\alpha}}(C)$, and by $H_{i+1}$ the open half space in $\mathbb{R}^n$ determined 
by $\mu_{\boldsymbol{\alpha}}(C_i) < \mu_{\boldsymbol{\alpha}}(C)$ for $ i = 0, ..., n-1$. Note that $h_i$ and $H_i$ are completely determined by the type of a chain.

\begin{rmk}
If a chain of type $\boldsymbol{t} = (r_0, r_1, ..., r_n; d_0, d_1, ..., d_n)$ is $\boldsymbol{\alpha}$-stable 
for some $\boldsymbol{\alpha} \in \mathbb{R}^n$, then $\boldsymbol{\alpha} \in H:= \bigcap_{i=1}^n H_i.$
\end{rmk}

The region $H$ is an $n$-dimensional open unbounded convex polyhedral cone with vertex
$\boldsymbol{v_0}:=(\frac{d_0}{r_0}-\frac{d_1}{r_1},...,\frac{d_0}{r_0}-\frac{d_n}{r_n})$. 
Indeed, the coefficients of the $\alpha_i$'s form the matrix
\[
\left(\begin{array}{ccccc}-r_1r_0&-r_2r_0&-r_3r_0&...&-r_nr_0 \\
r_1(r_2+...+r_n)&-r_2(r_0+r_1)&-r_3(r_0+r_1)&...&-r_n(r_0+r_1) \\
r_1(r_3+...+r_n)&r_2(r_3+...+r_n)&-r_3(r_0+...+r_2)&...&-r_n(r_0+...+r_2)\\
\vdots&\vdots&\vdots&...&\vdots\\
r_1r_n&r_2r_n&r_3r_n&...&-r_n(r_0+...+r_{n-1})\end{array}\right)
\]
whose determinant is $(-1)^nr_0...r_n(r_0+...+r_n)^{n-1}$. This is computed by elementary row operations.
Moreover, $\bigcap_{i=1}^n h_i = \{ \boldsymbol{v_0} \}$. We also have
%If $\boldsymbol{\alpha}$ and $\boldsymbol{\beta}$ are in $H$, then
\begin{flushleft}
i) the half line defined by $\boldsymbol{v_0} + t(\boldsymbol{\alpha} - \boldsymbol{v_0})$ with $t > 0$ is contained in $H$ for all $\boldsymbol{\alpha} \in H$, and \\
ii) the line segment defined by $(1-t) \boldsymbol{\alpha} + t \boldsymbol{\beta}$ with $ 0 \leq t \leq 1$ is contained in $H$ for all $\boldsymbol{\alpha}, \boldsymbol{\beta} \in H$.   
\end{flushleft}
Indeed, let
\begin{align*}
f_i(\boldsymbol\alpha)=f_i((\alpha_1,...,\alpha_n)) &=  \sum_{k=i+1}^nr_k  \sum_{k=1}^i \alpha_k r_k  - \sum_{k=0}^i r_k \sum_{k=i+1}^n \alpha_k r_k , \\
c_i &= \sum_{k=0}^i r_k \sum_{k=i+1}^n d_k  - \sum_{k=i+1}^n r_k \sum_{k=0}^i d_k.
\end{align*}
for $i = 0, 1, ..., n-1$.
Then $\boldsymbol{\alpha} \in H$ if and only if $f_i(\boldsymbol{\alpha}) < c_i$ and $f_i(\boldsymbol{v_0}) = c_i$ for all $i = 0, 1, ..., n-1$. Now we prove i) and ii) above.
Let $\boldsymbol{\alpha, \beta} \in H$. Then $f_i(\boldsymbol\alpha) < c_i$ and $f_i(\boldsymbol\beta) < c_i$. For each $i$,
i) $f_i(\boldsymbol{v_0}+t(\boldsymbol{\alpha}-\boldsymbol{v_0})) = (1-t)f_i(\boldsymbol{v_0}) + tf_i(\boldsymbol{\alpha}) = (1-t)c_i + t f_i(\boldsymbol{\alpha}) < (1-t)c_i + t c_i = c_i$, and ii)
$f_i((1-t) \boldsymbol{\alpha} + t \boldsymbol{\beta}) = (1-t)f_i(\boldsymbol{\alpha}) + t f_i(\boldsymbol{\beta}) < (1-t)c_i + t c_i = c_i$. Thus, $H$ is an open unbounded convex polyhedral cone in $\mathbb R^n$.

Note also that $H$ is nonempty because it contains the half line defined by $\boldsymbol{v_0} + t (0,...,0,1)$ for $t > 0$. 
Indeed $f_i(\boldsymbol{v_0} + t (0,...,0,1)) = f_i (\boldsymbol{v_0}) + t f_i((0,...,0,1)) = c_i - t r_n(r_0+...+r_i) < c_i $ for $i= 0,...,n-1$.

%rank one chains
\section{Stability of chains of type $\boldsymbol{t} = (1,...,1;d_0,...,d_n)$}

From now on we assume that $\boldsymbol{t} = (1,...,1;d_0,...,d_n)$ unless otherwise stated. Then the open convex polyhedral cone $H$ is determined by the inequalities
\[
(n - i)  \sum_{k=1}^i \alpha_k  - (i+1) \sum_{k=i+1}^n \alpha_k  < 
(i + 1 )\sum_{k=i+1}^n d_k  - (n-i) \sum_{k=0}^i d_k 
\]
for $ i = 0, ..., n-1$. 

\begin{prop} A chain $C = (L_0,...,L_n;\phi_1,...,\phi_n)$ of type ${\bf t}$ is $\boldsymbol{\alpha}$-stable for some $\boldsymbol{\alpha} \in \mathbb R^n$ if and only if $\phi_i \neq 0$ for all $i = 1,...,n$ and $\boldsymbol\alpha \in H$.
\end{prop}
\begin{proof} 
Suppose that $C$ is $\boldsymbol{\alpha}$-stable for some $\boldsymbol{\alpha} \in \mathbb{R}^n$. Then $\boldsymbol\alpha$ must be in $H$. Now we show $\phi_i\ne0$ for all $i=1, ..., n$. If $\phi_j = 0$ for some $j$, then the quotient chain $Q_j := C/C_j = (0,...,0,L_j,...,L_n;0,...0,\phi_{j+1},...,\phi_n)$ is a subchain, so $\mu_{\boldsymbol\alpha}(Q_j)<\mu_{\boldsymbol\alpha}(C)$, and since $Q_j$ is a quotient chain of an $\boldsymbol{\alpha}$-stable chain, we have $\mu_{\boldsymbol\alpha}(C)<\mu_{\boldsymbol\alpha}(Q_j)$. This is a contradiction.

Conversely suppose that $\boldsymbol\alpha\in H$ and $\phi_i \ne 0$ for all $i = 1,...,n$. Then
the $n$ standard subchains $C_0, C_1,...,C_{n-1}$ are the only subchains of $C$ composed of subbundles. Hence, $C$ is $\boldsymbol{\alpha}$-stable if and only if $\mu_{\boldsymbol{\alpha}}(C_i) < \mu_{\boldsymbol{\alpha}}(C)$ for all
$i = 0, 1, ...,n-1$, and $\mu_{\boldsymbol{\alpha}}(C_i) < \mu_{\boldsymbol{\alpha}}(C)$ for all $i = 0, 1, ...,n-1$ if and only if $\boldsymbol\alpha\in H$. Hence, $C$ is $\boldsymbol{\alpha}$-stable.

%both $C_j$ and $Q_j = (0,...,0,L_j,...,L_n;0,...0,\phi_{j+1},...,\phi_n)$ are subchains, so $\mu_{\boldsymbol{\alpha}}(C_j) < \mu_{\boldsymbol{\alpha}}(C)$ and 
%$\mu_{\boldsymbol{\alpha}}(Q_j) < \mu_{\boldsymbol{\alpha}}(C)$, i.e.
%\[
%(n - j) \sum_{k = 1}^j \alpha_k  - (j + 1) \sum_{k = j + 1}^n \alpha_k  < 
%(j + 1 )\sum_{k = j + 1}^n d_k  - (n-j) \sum_{k = 0}^j d_k , \tu{and}
%\]
%
%\[
%(n - j)  \sum_{k = 1}^j \alpha_k  - (j + 1)\sum_{k = j + 1}^n \alpha_k  > 
%(j + 1 )\sum_{k = j + 1}^n d_k  - (n-j) \sum_{k = 0}^j d_k ,
%\]
%respectively. This is a contradiction.
\end{proof}

\begin{prop}
Assume that $d_0 \geq d_1 \geq ... \geq d_n$. Then there exists an $\boldsymbol{\alpha}$-stable chain of type 
$\boldsymbol{t} = (1, 1, ..., 1;d_0, d_1, ..., d_n)$ for all $\boldsymbol{\alpha}$ in $H$.
\end{prop}
\begin{proof}
Since $d_{i-1} \geq d_i$, there exist line bundles $L_{i-1}, L_i$ of degrees $d_{i-1}, d_i$, respectively, and a nonzero map $\phi_i: L_i \rightarrow L_{i-1}$ for $i = 1, 2,...,n$.
The existence of the nonzero map is a well-known result of the Brill-Noether Theory. 
By Proposition 3.1, for all $\boldsymbol{\alpha}$ in $H$, the chain $(L_0,L_1,...,L_n;\phi_1,...,\phi_n)$ is $\boldsymbol{\alpha}$-stable of type $\boldsymbol t$.
\end{proof}

\begin{rmk} i) If $\boldsymbol{\alpha}$ is not in $H$, then there does not exist $\boldsymbol{\alpha}$-stable chain of type
$\boldsymbol{t}$. \\
ii) The region $H$ is called a chamber, where the $\boldsymbol{\alpha}$-stablility remains the same in the region.
\end{rmk}

\begin{prop}
Let $C$ be a chain of type $\boldsymbol{t} = (1, 1, ..., 1;d_0, d_1, ..., d_n)$.
If $C$ is $\boldsymbol{\alpha}$-semistable for some $\boldsymbol{\alpha}$ in $H$, then $C$ is  $\boldsymbol{\alpha}$-stable. In other words, on $H$, $C$ is  $\boldsymbol{\alpha}$-semistable if and only if $C$ is  $\boldsymbol{\alpha}$-stable.
\end{prop}
\begin{proof}
If $C$ is $\boldsymbol\alpha$-semistabe for some $\boldsymbol\alpha\in H$, then $\mu_{\boldsymbol\alpha}(C_i)\le\mu_{\boldsymbol\alpha}(C)$ for all $i=0,1,...,n-1$. Since $\boldsymbol\alpha\in H$, $\mu_{\boldsymbol\alpha}(C_i) < \mu_{\boldsymbol\alpha}(C)$ for all $i=0,1,...,n-1$. So, we are done. %, so we have $\mu_{\boldsymbol\alpha}(C)<\mu_{\boldsymbol\alpha}(C_i)$ for all $i=0,1,...,n-1$.

%Now let $C=(L_0,L_1,...,L_n;\phi_1,...,\phi_n)$. If $\phi_i\ne0$ for all $i=1,...,n$, then we are done.

%If $\phi_j = 0$ for some $j$, then both $C_j$ and $Q_j := C/C_j$ are subchains. Also note that $C_j = C/Q_j$. So, by semistability $\mu_{\boldsymbol{\alpha}}(C_j) \le \mu_{\boldsymbol{\alpha}}(C) \le \mu_{\boldsymbol{\alpha}}(C_j)$, so $\mu_{\boldsymbol{\alpha}}(C) = \mu_{\boldsymbol{\alpha}}(C_j)$. This implies that $\boldsymbol{\alpha}$ is not in $H$.
\end{proof}

\begin{defn} A chamber is a locally closed subset of  $\mathbb{R}^n$ where $\boldsymbol{\alpha}$-stability remains the same.
\end{defn}

\section{The moduli space of chains of type $\boldsymbol{t} = (1,...,1;d_0,...,d_n)$}
  
\subsection{The moduli spaces of $\boldsymbol{\alpha}$-stable chains}

Let $\mathcal{M}^s_{\boldsymbol{\alpha}}(\boldsymbol{t})$ be the moduli space of $\boldsymbol{\alpha}$-stable
chains of type $\boldsymbol{t} = (r_0,...,r_n;d_0,...,d_n)$ over a compact Riemann surface $X$ of genus $g$.
By the deformation theory of holomorphic chains
(\cite{AGS}, Theorem 3.8), the dimension of
$\mathcal{M}^s_{\boldsymbol{\alpha}}(\boldsymbol{t})$ at a smooth point $C$ is given by the formula
\[
\dim_C\mathcal{M}_{\boldsymbol{\alpha}}^s(\boldsymbol{t}) = (g-1)\left(\sum_{i=0}^n
r_i^2 - \sum_{i=1}^n r_ir_{i-1}\right) + \sum_{i=1}^n(r_id_{i-1} - r_{i-1}d_i) + 1.
\]
Moreover, if $\phi_i : E_i \rightarrow E_{i-1}$ is injective or generically surjective for all $i = 1,..., n$, then the
chain $C$ defines a smooth point of $\mathcal{M}^s_{\boldsymbol{\alpha}}({\bf t})$.

From now on we assume that $\boldsymbol{t} = (1,1,...,1;d_0,d_1,...,d_n)$. Then
$dim_C\mathcal{M}_{\boldsymbol{\alpha}}^s({\bf t}) = g + d_0 - d_n$; 
and if a chain $C = (L_0,L_1,...,L_n;\phi_1,...,\phi_n)$ 
of type $\boldsymbol{t}$ is $\boldsymbol{\alpha}$-stable, 
then by Proposition 2.1 $\phi_i \neq 0$ for all $i = 1,...,n$, so each homomorphism $\phi_i:L_i \rightarrow L_{i-1}$ is injective and it defines a smooth point in the moduli space $\mathcal{M}_{\boldsymbol{\alpha}}^s(\boldsymbol{t})$. We have the following result.
\begin{prop} If $\mathcal{M}^s_{\boldsymbol{\alpha}}(\boldsymbol{t}) \neq \emptyset$, then it is a smooth projective variety of dimension $g + d_0 - d_n$.
\end{prop}

\begin{lem} Let
 $C=(L_0,...,L_n;\phi_1,...,\phi_n)$ and $C'=(L_0,...,L_n;\psi_1,...,\psi_n)$ be two chains of
 type $\boldsymbol{t}$. Then
$C$ and $C'$  are isomorphic if and only if there exists $(\lambda_1,...,\lambda_n)\in(\mathbb{C}^*)^n$
such that $(\psi_1,...,\psi_n)=(\lambda_1\phi_1,...,\lambda_n\phi_n)$.
\end{lem}
\begin{proof} The automorphism group $Aut(L_i)$ is isomorphic to
$\mathbb{C}^*$ for $i=0,...,n$. Two chains $C$ and $C'$ are
isomorphic if and only if there exists a tuple $(t_0,...,t_n)\in(\mathbb{C}^*)^{n+1}$
such that $(\psi_1,...,\psi_n)=(t_0\phi_1 t_1^{-1},...,t_{n-1}\phi_n t_n^{-1})$.
Thus $(\psi_1,...,\psi_n)=(\lambda_1\phi_1,...,\lambda_n\phi_n)$
by letting $\lambda_i=t_{i-1}t_i^{-1}$ for $i=1,...,n$.
Conversely, if $(\psi_1,...,\psi_n)=(\lambda_1\phi_1,...,\lambda_n\phi_n)$, then
the $(n+1)$-tuple $(1,\lambda_1^{-1},\lambda_1^{-1}\lambda_2^{-1}...,\lambda_1^{-1}...\lambda_n^{-1})$
defines an isomorphism from $C$ to $C'$.
\end{proof}

Let $\mathcal L_i$ be the poincar\'e bundle over $Pic^{d_{i-1}-d_i}(X) \times X$. The pushforward $\nu^i_\ast \mathcal L_i$ is a coherent sheaf over $Pic^{d_{i-1}-d_i}(X)$, where $\nu^i$ is the projection map $Pic^{d_{i-1}-d_i}(X) \times X \to Pic^{d_{i-1}-d_i}(X)$. Now consider the pullback sheaf $\varphi_i^\ast\nu^i_\ast \mathcal L_i$, where $\varphi_i: Pic^{d_0}(X) \times Pic^{d_1}(X) \times \cdots \times Pic^{d_n}(X) \to Pic^{d_{i-1}-d_i}(X)$ is the map defined by $\varphi_i(L_0,...,L_n) = L_{i-1}L_i^{-1}$. %Let $Z=Pic^{d_0}(X) \times Pic^{d_1}(X) \times \cdots \times Pic^{d_n}(X)$.

\begin{thm}\label{stablemoduli} Let $\boldsymbol{t} = (1,1,...,1;d_0,d_1,...,d_n)$ and $X$ be a smooth curve of genus $g$. 
Assume that $d_{i-1} - d_i > 2g - 2$ for all $i = 1,...,n$. Then for $\boldsymbol{\alpha}$ in $H$
\[
\mathcal{M}^s_{\boldsymbol{\alpha}}(\boldsymbol{t}) 
= \mathbb P E_1\times_Z \mathbb P E_2 \times_Z \cdots \times_Z \mathbb P E_n,
%Pic^{d_0}(X) \times Pic^{d_1}(X) \times... \times Pic^{d_n}(X) \times \mathbb{P}^{d_0-d_1-g} \times ... \times \mathbb{P}^{d_{n-1}-d_n-g}
\]
where $E_i = \varphi_i^\ast\nu^i_\ast \mathcal L_i$ for $i=1,2,...,n$.
\end{thm}
\begin{proof}
Since $d_{i-1} - d_i > 2g - 2$, the direct image of a Poincar\'{e} line bundle over $X \times Pic^{d_{i-1}-d_i}(X)$ to $Pic^{d_{i-1}-d_i}(X)$ is a vector bundle over $Pic^{d_{i-1}-d_i}(X)$ for each $i = 1,...,n$. The fiber of the direct image of the Poincar\'{e} line bundle at $L \in Pic^{d_{i-1}-d_i}(X)$ is $H^0(X, L) \cong \mathbb C^{d_{i-1}-d_i + 1 - g}$ by the Riemann-Roch Theorem.
%For a point $L$ in $Pic^{d_{i-1}-d_i}(X)$, 
Then pulling back the direct image by the map $(L_0,L_1,...,L_n) \mapsto L_{i-1}L_i^{-1}$ from
\[
 Pic^{d_0}(X) \times Pic^{d_1}(X) \times \cdots \times Pic^{d_n}(X) \text{ to } Pic^{d_{i-1}-d_i}(X)
\]
is a vector bundle over $Pic^{d_0}(X) \times Pic^{d_1}(X) \times... \times Pic^{d_n}(X)$. Let $Z=Pic^{d_0}(X) \times Pic^{d_1}(X) \times \cdots \times Pic^{d_n}(X)$ and let $E_i$ be the pullback bundle over $Z$ for $i = 1,...,n$. Then ${E_1}\oplus \cdots \oplus E_n$ parameterizes all $\boldsymbol{\alpha}$-stable chains of type $\boldsymbol{t}$, and by Lemma 4.1 we have $\mathcal{M}^s_{\boldsymbol{\alpha}}(\boldsymbol{t}) 
= \mathbb{P}E_1 \times_Z \cdots \times_Z \mathbb{P}E_n$.

%Since $E_i$ a pullback of a direct image of a line bundle,
%\[
%\mathbb{P}E_i = Pic^{d_0}(X) \times Pic^{d_1}(X) \times... \times Pic^{d_n}(X) \times \mathbb{P}^{d_{i-1}-d_i-g}.
%\]
%Here $dim(H^0(X,L)) = d_{i-1}-d_i-g + 1$ by the Riemann-Roch Theorem. Now 
%Thus we have
%\[
%\mathcal{M}^s_{\boldsymbol{\alpha}}(\boldsymbol{t}) 
%= \mathbb{P}E_1 \times_Z ... \times_Z \mathbb{P}E_n = Z \times  \mathbb{P}^{d_0-d_1-g} \times ... \times \mathbb{P}^{d_{n-1}-d_n-g},
%\]
%where $Z = Pic^{d_0}(X) \times Pic^{d_1}(X) \times... \times Pic^{d_n}(X)$.

\end{proof}

If $d_{i-1} - d_i \leq 2g - 2$, then the direct image $\nu^i_\ast\mathcal L_i$ of a Poincar\'{e} line bundle over $X \times Pic^{d_{i-1}-d_i}(X)$ 
to $Pic^{d_{i-1}-d_i}(X)$ is not a vector bundle over $Pic^{d_{i-1}-d_i}(X)$ except in the case $g = 0$. 
However, if the direct image is restricted to a certain locally closed subset then it becomes
a vector bundle over the locally closed subset. For each $i$, denote by $U^{r_i}_{d_{i-1}-d_i}$ the locally closed subset of $Pic^{d_{i-1}-d_i}(X)$ 
whose points are line bundles with exactly $r_i+1$ sections, so that
\[
 U^{r_i}_{d_{i-1}-d_i} = W^{r_i}_{d_{i-1}-d_i} \setminus W^{r_i+1}_{d_{i-1}-d_i},
\]
where $W^r_d$ is the Brill-Noether locus of degree $d$, with at least $r + 1$ sections.

\begin{rmk}
 The locally closed subset $U^{r_i}_{d_{i-1}-d_i}$ is a dense open subset of the Brill-Noether locus $W^{r_i}_{d_{i-1}-d_i}$.
\end{rmk}

For each $i$, denote by $\pi_i$ the map $(L_0,L_1,...,L_n) \mapsto L_{i-1}L_i^{-1}$ from $ Z= Pic^{d_0}(X) \times Pic^{d_1}(X) \times \cdots \times Pic^{d_n}(X)$ to $Pic^{d_{i-1}-d_i}(X)$.

Now the direct image of a Poincar\'{e} line bundle over $X \times Pic^{d_{i-1}-d_i}(X)$ to $Pic^{d_{i-1}-d_i}(X)$   
is a vector bundle over $U^{r_i}_{d_{i-1}-d_i}$ for each $i$. Pulling back to $\pi_i^{-1}(U^{r_i}_{d_{i-1}-d_i})$ gives a vector
bundle over it. Denote by $E_{r_i}$ the vector bundles for $i=1,...,n$. Then $E_{r_1}\oplus\cdots \oplus E_{r_n}$ is a vector bundle over $Z_{r_1,...,r_n} := \bigcap_{i=1}^n \pi_i^{-1}(U^{r_i}_{d_{i-1}-d_i})$, and it parametrizes all holomorphic chains $(L_0,...,L_n;\phi_1,...,\phi_n)$ with $\dim (H^0(L_{i-1}L_i^{-1}))=r_i+1$ for all $i=1,...,n$.

\begin{rmk}
 The product $Z$ of the Jacobian varieties is the disjoint union of $Z_{r_1,...,r_n}$, i.e.,
\[
 Z = \underset{r_1 \geq 0,\cdots,r_n \geq 0}{\bigsqcup} Z_{r_1,...,r_n}.
\]
The component $Z_{r_1,...,r_n}$ consists of $(L_0,L_1,...,L_n)$ with $dim(H^0(L_i^{-1}L_{i-1})) = r_i + 1$ for all $i = 1,...,n$.
\end{rmk}

\begin{thm}
 Let $\boldsymbol{t} = (1,1,...,1;d_0,d_1,...,d_n)$ and $X$ be a smooth curve of genus $g$. Then 
for any $\boldsymbol{\alpha}$ in $H$, the moduli space 
$\mathcal{M}^s_{\boldsymbol{\alpha}}(\boldsymbol{t})$ is the disjoint union of the locally closed subsets $\mathbb P E_{r_1} \times_{Z_{r_1,...,r_n}} \cdots \times_{Z_{r_1,...,r_n}} \mathbb P E_{r_n}$,
i.e.,
\[
\mathcal{M}^s_{\boldsymbol{\alpha}}(\boldsymbol{t}) 
= \underset{r_1 \geq 0,...,r_n \geq 0}{\bigsqcup} \left( \mathbb P E_{r_1} \times_{Z_{r_1,...,r_n}} \cdots \times_{Z_{r_1,...,r_n}} \mathbb P E_{r_n} \right).
\]
\end{thm}
\begin{proof}
The fiber product
 $\mathbb{P}E_{r_1} \times_{Z_{r_1,...,r_n}} ... \times_{Z_{r_1,...,r_n}} \mathbb{P}E_{r_n}$ parameterizes the isomorphism classes of all chains with $\dim(H^0(L_i^{-1}L_{i-1})) = r_i + 1$ for all $i = 1,...,n$, and each $\boldsymbol\alpha$-stable chain has one of these properties. 
So, $\mathcal{M}^s_{\boldsymbol{\alpha}}(\boldsymbol{t})$ is the disjoint union of the fiber products over all $n$-tuples $(r_1,...,r_n)$ with $r_1 \geq 0,...,r_n \geq 0$.
\end{proof}

\begin{lem}
The Euler characteristic of $Z_{r_1,...,r_n}$ is $0$.
\end{lem}
\begin{proof}
 Consider the map $\pi : Z_{r_1,...,r_n} \rightarrow U^{r_1}_{d_0-d_1} \times \cdots \times U^{r_n}_{d_{n-1}-d_n}$ which is the restriction of $(\pi_1, ..., \pi_n): Z=Pic^{d_0}(X)\times\cdots\times Pic^{d_n}(X) \to Pic^{d_0-d_1}(X)\times\cdots\times Pic^{d_{n-1}-d_n}(X)$ to the locally open subset $Z_{r_1,...,r_n}$, i.e., $\pi(L_0,...,L_n)=(\pi_1(L_0,...,L_n),...,\pi_n(L_0,...,L_n))=(L_0L_1^{-1},..., L_{n-1}L_n^{-1})$. Then $\pi$ is onto and the fiber is isomorphic to $Pic^{d_n}(X)$. Given $(M_1,...,M_n) \in U^{r_1}_{d_{0}-d_1}\times \cdots \times U^{r_n}_{d_{n-1}-d_n}$ we have
\[
\pi(M_1...M_nL_n, M_2...M_nL_n, ...,M_nL_n, L_n) = (M_1,...,M_n) 
\]
for any $L_n \in Pic^{d_n}(X)$. So, $\chi(Z_{r_1,...,r_n}) = \chi(Pic^{d_n}(X))\chi(U^{r_1}_{d_{0}-d_1}\times \cdots \times U^{r_n}_{d_{n-1}-d_n}) = 0$.
 
%defined  by 
%\[ 
%\pi((L_0,...,L_n)) = \prod_{i=1}^n \pi_i((L_0,...,L_n)) = (L_0L_1^{-1},...,L_{n-1}L_n^{-1}). 
%\]

\end{proof}

\begin{thm}
 $\chi(\mathcal{M}^s_{\boldsymbol{\alpha}}(\boldsymbol{t})) = 0$.
\end{thm}
\begin{proof}
 By Theorem 4.2, $\chi(\mathcal{M}^s_{\boldsymbol{\alpha}}(\boldsymbol{t})) 
= \underset{r_1 \geq 0, \cdots ,r_n \geq 0}{\sum} \chi(\mathbb P E_{r_1} \times_{Z_{r_1,...,r_n}} \cdots \times_{Z_{r_1,...,r_n}} \mathbb P E_{r_n})$.
 
 Since $\mathbb P E_{r_1} \times_{Z_{r_1,...,r_n}} \cdots \times_{Z_{r_1,...,r_n}} \mathbb P E_{r_n}$ is a projective space bundle with base $Z_{r_1,...,r_n}$ and $\chi(Z_{r_1,...,r_n}) = 0$, then $\chi(\mathcal{M}^s_{\boldsymbol{\alpha}}(\boldsymbol{t})) = 0 $.
 
% $\chi(\mathcal{M}^s_{\boldsymbol{\alpha}}(\boldsymbol{t})) 
%= \underset{r_1 \geq 0,...,r_n \geq 0}{\sum} \chi(Z_{r_1,...,r_n}) \chi(\mathbb{P}^{r_1} \times ... \times \mathbb{P}^{r_n}) = 0$. 
\end{proof}

\subsection{The moduli spaces of $\boldsymbol{\alpha}$-semistable chains}

For chains of any type we have the following lemma.
\begin{lem}\label{barycenter}
 If $0 \rightarrow C' \rightarrow C \rightarrow C'' \rightarrow 0$ is an exact sequence of chains, then given $\boldsymbol{\alpha}$, 
$\mu_{\boldsymbol{\alpha}}(C)$ is the barycenter of $\mu_{\boldsymbol{\alpha}}(C')$ and $\mu_{\boldsymbol{\alpha}}(C'')$.
\end{lem}
\begin{proof}
 Let $\boldsymbol{t} = (d_0,...,d_n;r_0,...,r_n)$, $\boldsymbol{t'} = (d_0',...,d_n';r_0',...,r_n')$, 
$\boldsymbol{t''} = (d_0'',...,d_n'';r_0'',...,r_n'')$ be the types of $C$, $C'$, $C''$, respectively. Let $d = \sum_{i=0}^n d_i$, 
$d'= \sum_{i=0}^n d_i'$, $d''= \sum_{i=0}^n d_i''$, $r = \sum_{i=0}^n r_i$, $r'= \sum_{i=0}^n r_i'$, $r''= \sum_{i=0}^n r_i''$.
Then $d = d' + d''$ and $r = r' + r''$. Given $\boldsymbol{\alpha} = (\alpha_1,...,\alpha_n)$,
\begin{align*}
 \mu_{\boldsymbol{\alpha}}(C) &= \frac{d+\sum_{i=1}^n\alpha_ir_i}{r} = \frac{(d'+d'')+\sum_{i=1}^n\alpha_i(r_i'+r_i'')}{r' + r''} \\
&= \frac{(d'+ \sum_{i=1}^n\alpha_i r_i') + (d'' + \sum_{i=1}^n\alpha_i r_i'' ) }{r' + r''} 
= \frac{ r'\frac{(d'+ \sum_{i=1}^n\alpha_i r_i')}{r'} + r''\frac{(d'' + \sum_{i=1}^n\alpha_i r_i'')}{r''}}{r' + r''}\\
&= \frac{r'\mu_{\boldsymbol{\alpha}}(C') + r''\mu_{\boldsymbol{\alpha}}(C'')}{r' + r''}.
\end{align*}
\end{proof}

Recall that for a chain $C$, $h_{i+1}=\{ \boldsymbol\alpha : \mu_{\boldsymbol\alpha}(C_i) = \mu_{\boldsymbol\alpha}(C)\}, \ H_{i+1}=\{ \boldsymbol\alpha : \mu_{\boldsymbol\alpha}(C_i) < \mu_{\boldsymbol\alpha}(C)\}$ for $i=0,1,...,n-1$, and $H=\bigcap_{i=1}^n H_i$.

From now on we assume that chains are of type $\boldsymbol t =(1,...,1;d_0,...,d_n)$. If a chain $C$ is $\boldsymbol{\alpha}$-stable for some $\boldsymbol{\alpha}$ in $H$, then $C$ is $\boldsymbol{\alpha}$-semistable but not $\boldsymbol{\alpha}$-stable for $\boldsymbol{\alpha}\in \overline{H} \setminus H$.

For $I \subseteq \{1,  ... ,n \}$, define
\[
\mathcal C_I := \left( \bigcap_{i \in I} h_i \setminus \bigcup_{|J|=|I|+1} \left(\bigcap_{ j \in J } h_j  \right) \right) \bigcap \overline{H}. 
\]
The subset $\mathcal C_I$ is an $|I|$-dimensional locally closed subset of $\overline{H}$. 
If $I = \{1,...,n\}$, then we define $\mathcal C_I := \{ \boldsymbol{v_0} \}$, and if $I = \emptyset$, then we define $\mathcal C_I := H$.

\begin{lem}\label{semi} The chain $C=(L_0,...,L_n;\phi_1,...,\phi_n)$ with $\phi_i \neq 0$ for all $i$ is $\boldsymbol{\alpha}$-stable for $\boldsymbol{\alpha} \in H$. If $\boldsymbol{\alpha} \in \mathcal C_{\{i_0\}}$, then it is $\boldsymbol{\alpha}$-semistable and $C$ is S-equivalent to the chain $C'=(L_0,...,L_n;\psi_1,...,\psi_n)$ with $\psi_i = \phi_i$ for all $i \neq i_0$ and $\psi_{i_0} = 0$.
\end{lem}

\begin{proof}
There exist Jordan-H\"{o}lder filterations $0 \subset C_{i_0-1} \subset C$ and $0 \subset C_{i_0-1} \subset C'$ where $C_{i_0-1}$ is the $(i_0-1)$-th standard subchain $(L_0,...,L_{i_0-1},0,...,0;\phi_1,...,\phi_{i_0-1},0,...,0)$. First, we show that $C_{i_0-1}$ is $\boldsymbol{\alpha}$-stable with slope $\mu_{\boldsymbol{\alpha}}(C)$ for $\boldsymbol{\alpha} \in \mathcal C_{\{i_0\}}$. It is enough to check $\boldsymbol{\alpha}$-stability against the standrad subchains $C_0, C_1, ..., C_{i_0-2}$. Since $\boldsymbol{\alpha}$ is in $h_{i_0}$ but not in $h_1, h_2, ...,h_{i_0-1}$, $\boldsymbol{\alpha}\in\bigcap_{i=0}^{i_0-2} H_i$ and so we have $\mu_{\boldsymbol{\alpha}}(C_i)<\mu_{\boldsymbol{\alpha}}(C)=\mu_{\boldsymbol{\alpha}}(C_{i_0-1})$ for $i=0,1,...,i_0-2$.

Second, we show that the quotient chains $C/C_{i_0-1} = C'/C_{i_0-1}$ are also $\boldsymbol{\alpha}$-stable with slope $\mu_{\boldsymbol{\alpha}}(C)$. These quotient chains are isomorphic to $(0,...,0,L_{i_0},...,L_n;0,...,0,\phi_{i_0+1},...,\phi_n)$.
Since $0 \rightarrow C_{i_0-1} \rightarrow C \rightarrow C/C_{i_0-1} \rightarrow 0$, 
by Lemma \ref{barycenter}, $\mu_{\boldsymbol{\alpha}}(C)$ is the barycenter of $\mu_{\boldsymbol{\alpha}}(C_{i_0-1})$ and $\mu_{\boldsymbol{\alpha}}(C/C_{i_0-1})$, but $\mu_{\boldsymbol{\alpha}}(C_{i_0-1})=\mu_{\boldsymbol{\alpha}}(C)$, so 
$\mu_{\boldsymbol{\alpha}}(C/C_{i_0-1}) = \mu_{\boldsymbol{\alpha}}(C)$. Now we show that they are $\boldsymbol{\alpha}$-stable. Any nonzero subchain of $C/C_{i_0}$ consiting of subbundles is in the form $C_j/C_{i_0}$ for some $j > i_0$, 
so we have an exact sequence
$0 \rightarrow C_{i_0} \rightarrow C_j \rightarrow C_j/C_{i_0} \rightarrow 0$. Since  
$\mu_{\boldsymbol{\alpha}}(C_j) < \mu_{\boldsymbol{\alpha}}(C) = \mu_{\boldsymbol{\alpha}}(C_{i_0})$,
by Lemma \ref{barycenter} $\mu_{\boldsymbol{\alpha}}(C_j/C_{i_0}) < \mu_{\boldsymbol{\alpha}}(C_j)$ and so, $\mu_{\boldsymbol{\alpha}}(C_j/C_{i_0}) < \mu_{\boldsymbol{\alpha}}(C_j) < \mu_{\boldsymbol{\alpha}}(C) 
= \mu_{\boldsymbol{\alpha}}(C/C_{i_0})$. Thus $C/C_{i_0}$ is $\boldsymbol{\alpha}$-stable.

Consequently, the $C$ and $C'$ are $S$-equivalent.

% is also $\boldsymbol{\alpha}$-stable with slope $\boldsymbol{\alpha}$-slope $\mu_{\boldsymbol{\alpha}}(C/C_{i_0}) = \mu_{\boldsymbol{\alpha}}(C'/C'_{i_0}) = \mu_{\boldsymbol{\alpha}}(C)$. 
% Since $0 \rightarrow C_{i_0} \rightarrow C \rightarrow C/C_{i_0} \rightarrow 0$, by Lemma \ref{barycenter}, $\mu_{\boldsymbol{\alpha}}(C/C_{i_0}) = \mu_{\boldsymbol{\alpha}}(C)$.

% Any nonzero subchain of $C/C_{i_0}$ is in the form $C_j/C_{i_0}$ for some $j > i_0$, so we have an exact sequence $0 \rightarrow C_{i_0} \rightarrow C_j \rightarrow C_j/C_{i_0} \rightarrow 0$. Since $\mu_{\boldsymbol{\alpha}}(C_j) < \mu_{\boldsymbol{\alpha}}(C) = \mu_{\boldsymbol{\alpha}}(C_{i_0})$, by Lemma \ref{barycenter} $\mu_{\boldsymbol{\alpha}}(C_j/C_{i_0}) < \mu_{\boldsymbol{\alpha}}(C_j)$.
%So $\mu_{\boldsymbol{\alpha}}(C_j/C_{i_0}) < \mu_{\boldsymbol{\alpha}}(C_j) < \mu_{\boldsymbol{\alpha}}(C) = \mu_{\boldsymbol{\alpha}}(C/C_{i_0})$, and $C/C_{i_0}$ is $\boldsymbol{\alpha}$-stable.
%Similarly, $C'/C'_{i_0}$ is $\boldsymbol{\alpha}$-stable with $\mu_{\boldsymbol{\alpha}}(C'/C'_{i_0}) = \mu_{\boldsymbol{\alpha}}(C') = \mu_{\boldsymbol{\alpha}}(C)$.
%We conclude that gr$(C) = (C_{i_0}/0)\oplus (C/C_{i_0})$ and gr$(C') = (C'_{i_0}/0)\oplus (C'/C_{i_0}')$, and they are isomorphic.
\end{proof}

\begin{lem}\label{semi2} If $\boldsymbol{\alpha} \in \mathcal C_I$, then $C$ is $\boldsymbol{\alpha}$-semistable and $C$ is S-equivalent to the chain $C'=(L_0,...,L_n;\psi_1,...,\psi_n)$ with $\psi_i = \phi_i$ for $i \neq I$ and $\psi_j = 0$ for $j \in I$.
\end{lem}
\begin{proof}
Let $I=\{i_1, ..., i_m\}$. Then their Jordan-H\"{o}lder filterations are $0\subset C_{i_1-1} \subset \cdots \subset C_{i_m-1} \subset C$ and $0\subset C_{i_1-1}' \subset \cdots \subset C_{i_m-1}' \subset C'$ and their graduations gr$(C)=(C_{i_1-1}/0)\oplus  (C_{i_2-1}/C_{i_1-1})\oplus\cdots \oplus (C/C_{i_m-1})$ and gr$(C')=(C_{i_1-1}'/0)\oplus  (C_{i_2-1}'/C_{i_1-1}') \oplus \cdots \oplus(C'/C_{i_m-1}')$ are isomorphic.
\end{proof}

Denote by $\mathcal{M}_{\boldsymbol{\alpha}}(\boldsymbol{t})$ the moduli space of the $S$-equivalent classes of $\boldsymbol{\alpha}$-semistable chains of type $\boldsymbol{t}$.

%If $\boldsymbol{\alpha}$ moves to the boundary $\overline{H} \setminus H$, then $C$ is not $\boldsymbol{\alpha}$-stable but $\boldsymbol{\alpha}$-semistable.

\begin{lem}\label{semimoduli} Suppose $d_{i-1}-d_i>2g-2$ for all $i=1,...,n$.
 If $\boldsymbol{\alpha}$ is in $\mathcal C_{\{i_0\}}$, then
\[
\mathcal{M}_{\boldsymbol{\alpha}}(\boldsymbol{t}) 
= \mathbb P E_1 \times_Z \cdots \times_Z \widehat{\mathbb P E_{i_0}} \times_Z \cdots \times_Z \mathbb P E_n.
\]
\end{lem}
\begin{proof}
An $\boldsymbol \alpha$-semistabe chain $(L_0, L_1, ..., L_n; \phi_1, ... ,\phi_n)$ representing a fiber of the vector bundle $E_{i_0}$ at $(L_0, L_1, ..., L_n)$ constructed in the proof of Theorem \ref{stablemoduli} is $S$-equivalent 
to the chain
\[
(L_0, L_1, ..., L_n; \phi_1, ..., \phi_{i_0-1}, 0, \phi_{i_0+1}, ...,\phi_n).
\]
This means that each fiber in the component $\mathbb{P}E_{i_0}$ squashes to its base point in the moduli space $\mathcal{M}_{\boldsymbol{\alpha}}(\boldsymbol{t})$, so that the component $\mathbb{P}E_{i_0}$ becomes $Z=Pic^{d_0}(X) \times \cdots \times Pic^{d_n}(X)$. Thus 
\[
\mathcal{M}_{\boldsymbol{\alpha}}(\boldsymbol{t}) 
= \mathbb P E_1 \times_Z \cdots \times_Z Z \times_Z \cdots \times_Z \mathbb P E_n= \mathbb P E_1 \times_Z \cdots \times_Z \widehat{\mathbb P E_{i_0}} \times_Z \cdots \times_Z \mathbb P E_n.
\]
\end{proof}

%\begin{lem}\label{semimoduli}
% If $\boldsymbol{\alpha}$ is in $h_{i_0+1} \cap \overline{H}$, then
%\[
%\mathcal{M}_{\boldsymbol{\alpha}}(\boldsymbol{t}) 
%= \prod_{i=0}^n Pic^{d_i}(X) \times \prod_{i \neq i_0} \mathbb{P}^{d_{i-1}-d_i-g}.
%\]
%\end{lem}
%\begin{proof}
% The chain 
% \[ 
% (L_0, L_1, ..., L_n; \phi_1, ... ,\phi_n) 
% \] representing each fiber of the vector bundle $E_{i_0}$ constructed in the proof of Theorem \ref{stablemoduli} is $S$-equivalent 
%to the chain 
%\[
%(L_0, L_1, ..., L_n; \phi_1, ..., \phi_{i_0-1}, 0, \phi_{i_0+1}, ...,\phi_n).
%\]
%This means that each fiber squashes to its base point in the moduli space $\mathcal{M}_{\boldsymbol{\alpha}}(\boldsymbol{t})$, so that
%$\mathbb{P}E_{i_0} = \prod_{i=0}^n Pic^{d_i}(X)$. Thus 
%\[
% \mathbb{P}E_1 \times_Z ... \times_Z \mathbb{P}E_n = Z \times \prod_{i \neq i_0} \mathbb{P}^{d_{i-1}-d_i-g},
%\]
%where $Z=\prod_{i=0}^n Pic^{d_i}(X)$.
%\end{proof}

%For $I \subseteq \{1,  ... ,n \}$, define
%\[
%C_I := \left( \bigcap_{i \in I} h_i \setminus \bigcup_{|J|=|I|+1} \left(\bigcap_{ j \in J } h_j  \right) \right) \bigcap \overline{H}. 
%\]
%The subset $C_I$ is an $|I|$-dimensional locally closed subset of $\overline{H}$. 
%If $I = \{1,...,n\}$, then we define $C_I := \{ \boldsymbol{v_0} \}$, and if $I = \emptyset$, then we define $C_I := H$.

\begin{thm} Suppose $d_{i-1}-d_i>2g-2$ for all $i=1,...,n$. For $\boldsymbol{\alpha}$ in $C_I$, the moduli space $\mathcal{M}_{\boldsymbol{\alpha}}(\boldsymbol{t})$ of $\boldsymbol{\alpha}$-semistabe chains is the fiber product of $\mathbb P E_i$'s over $Z$ skipping the components $\mathbb P E_j$'s for $j \in I$. In particular, if $I=\emptyset$, then $\mathcal{M}_{\boldsymbol{\alpha}}(\boldsymbol{t})=\mathcal{M}^s_{\boldsymbol{\alpha}}(\boldsymbol{t})$ since $C_I=H$, and if $I=\{1,2,...,n\}$, then $\mathcal{M}_{\boldsymbol{\alpha}}(\boldsymbol{t})=pt$ since $C_I=\{\boldsymbol{v_0} \}$, and every chain is $S$-equivalent to $(L_0,L_1,...,L_n;0,...,0)$.
\end{thm}
\begin{proof}
This is by Lemma \ref{semi2} and \ref{semimoduli}.
\end{proof}

\begin{rmk} If $d_{i-1}-d_i < 2g-2$ for some $i$, then for $\boldsymbol{\alpha} \in C_I$, the moduli space $\mathcal{M}_{\boldsymbol{\alpha}}(\boldsymbol{t})$ of $\boldsymbol{\alpha}$-semistabe chains is the disjoint unit of the fiber products of $\mathbb P E_{r_i}$'s over $Z_{r_1...r_n}$ skipping the components $\mathbb PE_{r_j}$'s for $j \in I$.

\end{rmk}

\section{A variation of GIT and the parameter $\boldsymbol{\alpha}$}
Let $\boldsymbol{t}=(1,...,1;d_0,d_1,...,d_n)$, and we assume that $d_{i-1} - d_i > 2g - 2$ for $i = 1,...,n$. Then the moduli space $\mathcal M_\alpha^s(\boldsymbol{t})$ is $\mathbb P E_1 \times_Z \cdots \times_Z \mathbb P E_n$. Each $\mathbb P E_i$ is the projective space bundle coustructed in \cite{Har}. It is locally a projective space over $U\subset Z$. Here $U=\text{Spec}(A)$ for some ring $A$. In other words, $\mathbb P E_i|_U = \mathbb P_U^{r_i}$ where $r_i=d_{i-1}-d_i-g$. Then $\mathbb P_U^{r_i}=\text{Proj}\bigoplus_{d\ge 0} S_{d}$ where $S=A[x_0,x_1,...,x_{r_i}]$. We can reconstruct $\mathbb PE_i$ by gluing $\mathbb P_U^{r_i}$'s.

Let $G=G_m=\mathbb C^*$ and let $\chi_i: G\to \mathbb C^*$ be the character defined by $\chi_i(\lambda)=\lambda^{k_i}$ for some integer $k_i$. Let $S[E_i]^{G,\chi_i^m}$ be the set of polynomials over $A$ satisfying $f(g\cdot v)=\chi_i(g)^m f(v)$ for all $g \in G$ and $v\in \mathbb C^{r_i+1}$. Here the action $g\cdot v$ is the scalar multiplication of $\mathbb C^\ast$ on $\mathbb C^{r_i+1}$.
If $k_i>0$, then we have $\mathbb P_U^{r_i}\cong \text{Proj}\bigoplus_{m\ge0} S[E_i]^{G,\chi_i^m}$. If $k_i<0$, then $\text{Proj}\bigoplus_{m\ge0} S[E_i]^{G,\chi_i^m}=\text{Proj}(A) =\emptyset$, and if $k_i=0$, then $\text{Proj}\bigoplus_{m\ge0} S[E_i]^{G,\chi_i^m}=\text{Proj} \bigoplus_{m\ge0} A = \text{Proj}(A[t]) = \text{Spec}(A)=U$. We write $(U \times \mathbb C^{r_i+1}) /_{\chi_i} G := \text{Proj}\bigoplus_{m\ge0} S[E_i]^{G,\chi_i^m}$, and write $E_i/_{\chi_i} G$ for their gluing. If $k_i>0$ for $i=1,...,n$, then $E_i/_{\chi_i} G \cong \mathbb P E_i$ for $i=1,...,n$, and 
\[
\mathcal M_\alpha(\boldsymbol{t})=\mathcal M_\alpha^s(\boldsymbol{t})=\mathbb P E_1 \times_Z \cdots \times_Z \mathbb P E_n \cong (E_1/_{\chi_1} G) \times_Z \cdots \times_Z (E_n/_{\chi_n} G).
\]

Define $M =\{(k_1,...,k_n) \in \mathbb Z^n ~ | ~ k_1>0, ..., k_n>0\}$. Here $\mathbb Z \cong \text{Hom}(G,\mathbb G_m)$. Then there is a correspondence between $M'=\{(k_1,...,k_n) \in \mathbb Z^n ~ | ~ k_1\ge 0, ..., k_n\ge0\}$ and $\overline H$ which define the variation of the projective quotient and the variation of $\boldsymbol \alpha$, respectively.

For $I\subseteq \{1,2,...,n\}$ define $\mathcal D_I=\{(k_1,...,k_n)\in M' ~|~ k_i=0 \text{ for } i\in I \text{ and } k_i>0 \text{ for } i\notin I \}$. If $I=\{1,2,...,n\}$, then $\mathcal D_I=\{(0,0,...,0)\}$ and if $I=\emptyset$, then $\mathcal D_I=M$. The subsets $\mathcal D_I$ corresponds to the regions $\mathcal C_I$.

\begin{thm}
The parameter $\boldsymbol{\alpha}$ in $\mathbb R^n$ corresponds to the charaters in $\mathbb Z^n \cong  \text{Hom}(G,G_m)^n$.
\end{thm}
\begin{proof} If $\boldsymbol{\alpha} \notin \overline{H}$, then $\mathcal M_{\boldsymbol{\alpha}}(\boldsymbol{t})=\emptyset$, and if $(k_1,...,k_n)\notin M'$, then  $E_i/_{\chi_i} G=\emptyset$ for some $i$, so $(E_1/_{\chi_1} G) \times_Z \cdots \times_Z (E_n/_{\chi_n} G)=\emptyset$. If $\boldsymbol{\alpha} \in \mathcal C_I$, then the component $\mathbb P E_i$ collapses to $Z$ for $i\in I$ by the $S$-equivalence, so $\mathcal M_{\boldsymbol{\alpha}}(\boldsymbol{t})$ is the fiber product of $\mathbb P E_j$'s over $Z$ skipping the components $\mathbb P E_i$ for $i \in I$, and if $(k_1,...,k_n)\in \mathcal D_I$, then $E_i/_{\chi_i} G=Z$ for $i\in I$, and so $(E_1/_{\chi_1} G) \times_Z \cdots \times_Z (E_n/_{\chi_n} G)$ is the fiber product of $E_j/_{\chi_j} G$'s over $Z$ skipping the components $(E_i/_{\chi_i} G)$ for $i\in I$. Here $E_j/_{\chi_j} G\cong\mathbb P E_j$ for $j \notin I$.

Finally, if $I=\emptyset$, then $\mathcal D_I = M$ and $\mathcal C_I = H$. Then for $\boldsymbol \alpha \in H$, $\mathcal M_{\boldsymbol{\alpha}}(\boldsymbol{t})=\mathbb P E_1 \times_Z \cdots \times_Z \mathbb P E_n$, and for $(k_1,...,k_n) \in M$, $(E_1/_{\chi_1} G) \times_Z \cdots \times_Z (E_n/_{\chi_n} G)\cong \mathbb P E_1 \times_Z \cdots \times_Z \mathbb P E_n=\mathcal M_{\boldsymbol{\alpha}}(\boldsymbol{t})$. Thus $\mathcal D_I$ corresponds to $\mathcal C_I$ and they both define the same moduli space.
\end{proof}

\vskip 4cm

\section{Final remark}
If $d_{i-1} - d_i > 2g - 2$ for all $i$, then the pullback sheaf $\varphi_i^\ast\nu^i_\ast \mathcal L_i$ is locally free. So, Proj$(\varphi_i^\ast\nu^i_\ast \mathcal L_i)$ is the projective space bundle of the vector bundle. Then the fiber product of these over $Z$ can be identified with the moduli space $\mathcal M^s_{\boldsymbol{\alpha}}(\boldsymbol{t})$.

If $d_{j-1} - d_j < 2g - 2$ for some $j$, then the pullback sheaf $\varphi_j^\ast\nu^j_\ast \mathcal L_j$ is not locally free. The proj Proj$(\varphi_j^\ast\nu^j_\ast \mathcal L_j)$ of the sheaf $\varphi_j^\ast\nu^j_\ast \mathcal L_j$ parametrizes the isomorphic classes of holomorphic chains over locally closed subsets of $Z$.

\end{document}